\documentclass{amsart}
\usepackage{amsmath,amsfonts,amssymb,epsfig,graphicx,subfigure,graphs,amsopn,verbatim,color}

\newtheorem{rems}{Remark}

\newtheorem*{thm}{Theorem}

\newtheorem{lem}{Lemma}
\newtheorem*{prop}{Proposition} 

\theoremstyle{remark}
\newtheorem{example}{Example}
\begin{document}

\title[Random permutations and ergodicity for the Euler adic]{Random permutations and unique fully
supported ergodicity for the Euler adic transformation}
\author{Sarah Bailey Frick}
\email{frick@math.ohio-state.edu}
\author{Karl Petersen}
\email{petersen@math.unc.edu}
\address{Department of Mathematics,
         CB 3250, Phillips Hall,\\
         University of North Carolina,
         Chapel Hill, NC 27599 USA}


\begin{abstract}
There is only one fully supported ergodic invariant probability
measure for the adic transformation on the space of infinite paths
in the graph that underlies the Eulerian numbers.  This result may
partially justify a frequent assumption about the equidistribution
of random permutations.
\end{abstract}
\maketitle

\section{Introduction}
  We give a new proof of the main result of \cite{BKPS}, ergodicity of the symmetric measure for the Euler adic transformation. In fact
  we prove a stronger result by a different method: the symmetric measure is the {\em only} fully supported ergodic invariant probability measure for the Euler adic system.

This result has a probabilistic interpretation. For each
  $n=0,1,2,\dots,k=0,1,\dots n$, the
 \emph{Eulerian number} $A(n,k)$
 is the number of permutations
$i_1i_2\dots i_{n+1}$ of $\{1,\dots,n+1\}$ with exactly $k$
\emph{rises} (indices $j=1,2,\dots,n$ with $i_j<i_{j+1}$) and $n-k$
\emph{falls} (indices $j=1,2,\dots, n$ with $i_j>i_{j+1}$). (See,
for example, \cite{Comtet} for basic facts about the Eulerian
numbers.) Besides their combinatorial importance, these numbers are
also of interest in connection with the statistics of rankings: see,
for example, \cite{CKSS},\cite{Esseen},\cite{FLW},\cite{FL}, and
\cite{OV}.  In studying random permutations, it is often assumed
that all permutations are equally likely, each permutation of length
$n+1$ occurring with probability $1/(n+1)!$. Our main result implies
that in a sense (see Remark \ref{randomperms}) there is a unique way
to choose permutations at random so that (1) consistency holds: the
distribution on permutations of $\{1,\dots,n+m\}$ induces, upon
deletion of $\{n+1,\dots,m\}$, the distribution on $\{1,\dots,n\}$;
(2) any two permutations of the same length which have the same
number of rises are equally likely; and (3) every permutation has
positive probability.  This unique way is to make all permutations
of the same length, \emph{no matter the number of rises}, equally
likely.

This conclusion is achieved by means of ergodic theory, studying
invariant measures for the adic (Bratteli-Vershik) transformation on
the space of infinite paths on the Euler graph. We do not see how to
prove this result by a straightforward combinatorial argument,
because it is very difficult to make comparative asymptotic
estimates of Eulerian numbers. Indeed, the key to our argument is to
finesse this difficulty (in the proof of Proposition \ref{Unique})
by setting up a situation in which the dominant terms of the two
expressions being compared are identical.

Not having the space here to survey all the literature on adic
systems, we remark just that in certain senses they represent all
dynamical systems, also have strong connections with group
representations, $C^*$ algebras, probability theory, and
combinatorics, and thus facilitate the study of important problems
in those areas. See \cite{KV,GPS,HPS,PS,M,MP} and the references
that these works contain for background. Basic attributes of adic
systems are the various kinds of dimension groups, which provide
invariants, and the states or traces, which correspond to
adic-invariant, or ``central", measures---see, for example,
\cite{KV}. These have been determined for many stationary and
essentially simple adic systems and some classes of nonstationary
and nonsimple ones \cite{SBF1,SBF2}, but not yet for many key
individual examples. For the Pascal adic, the identification of the
ergodic invariant measures as the one-parameter family of Bernoulli
measures on the 2-shift is closely connected with de Finetti's
Theorem, the Hewitt-Savage 0,1 Law, and uniform distribution of the
Kakutani interval splitting procedure---see the discussion in
\cite{PS}. Identifying the ergodic measures for new natural examples
such as the Euler adic will also have consequences in other areas:
mentioned below are connections with random permutations (Remark
\ref{randomperms}) and reinforced random walks (Remark \ref{rw}).

\section{The Euler adic system}

 The \emph{Euler graph} $\Gamma$ is a Bratteli diagram with levels
 $n=0,1,2,\dots$.
There are $n+1$ vertices of each level $n$, labelled $(n,k)$, $0\leq
k\leq n$,
 and each vertex $(n,k)$ is connected to vertex $(n+1,k)$ by $k+1$ edges and
 to
vertex $(n+1,k+1)$ by $n-k+1$ edges.  The number of paths into any
vertex $(n,k)$ is the Eulerian number $A(n,k)$.  The Eulerian
numbers have the recursion

\begin{equation}
A(n,k)=(n-k+1)A(n-1,k-1)+(k+1)A(n-1,k), \label{recursion}
\end{equation}
where $A(0,0)=1$ and by convention $A(n,k)=0$ for $k\notin
\{0,1,...,n\}$.

\begin{figure}[h]
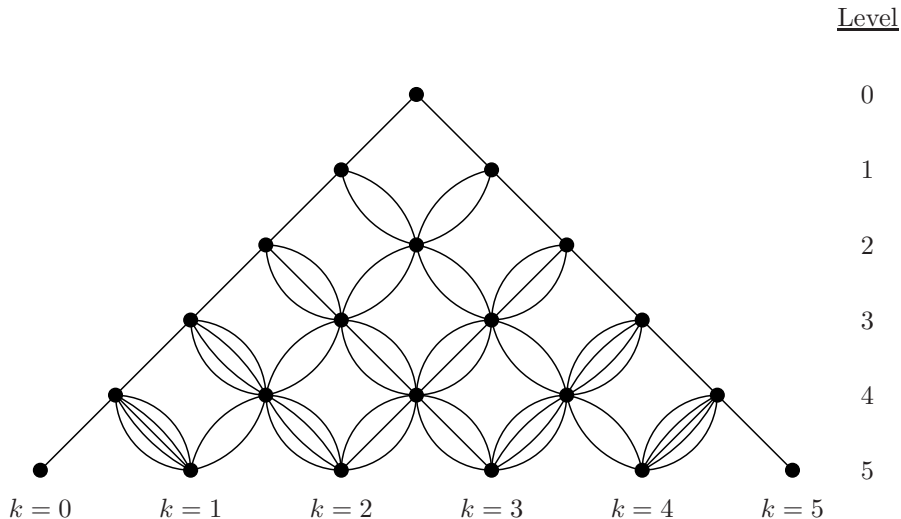

\begin{graph}(9,7)(0,-1.5) \roundnode{v00}(4,4)
\roundnode{v10}(3,3)\roundnode{v11}(5,3)
\roundnode{v20}(2,2)\roundnode{v21}(4,2)\roundnode{v22}(6,2)\edge{v00}{v10}\edge{v00}{v11}\edge{v10}{v20}\bow{v10}{v21}{.15}
\bow{v10}{v21}{-.15}\bow{v11}{v21}{-.15}
\bow{v11}{v21}{.15}\edge{v11}{v22}\roundnode{v30}(1,1)\roundnode{v31}(3,1)\roundnode{v32}(5,1)\roundnode{v33}(7,1)
\edge{v20}{v30}\edge{v20}{v31}\bow{v20}{v31}{.2}\bow{v20}{v31}{-.2}\bow{v21}{v31}{.15}\bow{v21}{v31}{-.15}
\bow{v21}{v32}{.15}\bow{v21}{v32}{-.15}\edge{v22}{v32}\bow{v22}{v32}{.2}\bow{v22}{v32}{-.2}\edge{v22}{v33}
\roundnode{v40}(0,0)\roundnode{v41}(2,0)\roundnode{v42}(4,0)\roundnode{v43}(6,0)\roundnode{v44}(8,0)
\edge{v30}{v40}\bow{v30}{v41}{.075}\bow{v30}{v41}{.2}\bow{v30}{v41}{-.075}\bow{v30}{v41}{-.2}\bow{v31}{v41}{.15}
\bow{v31}{v41}{-.15}\bow{v31}{v42}{.2}\bow{v31}{v42}{-.2}\edge{v31}{v42}\edge{v32}{v42}\bow{v32}{v42}{.2}\bow{v32}{v42}{-.2}
\bow{v32}{v43}{.15}\bow{v32}{v43}{-.15}
\bow{v33}{v43}{.075}\bow{v33}{v43}{.2}\bow{v33}{v43}{-.075}\bow{v33}{v43}{-.2}\edge{v33}{v44}
\roundnode{v50}(-1,-1)\roundnode{v51}(1,-1)\roundnode{v52}(3,-1)\roundnode{v53}(5,-1)\roundnode{v54}(7,-1)\roundnode{v55}(9,-1)
\edge{v40}{v50}\edge{v40}{v51}\bow{v40}{v51}{.075}\bow{v40}{v51}{.2}\bow{v40}{v51}{-.075}\bow{v40}{v51}{-.2}
\bow{v41}{v51}{.15}\bow{v41}{v51}{-.15}
\bow{v41}{v52}{.075}\bow{v41}{v52}{.2}\bow{v41}{v52}{-.075}\bow{v41}{v52}{-.2}
\edge{v42}{v52}\bow{v42}{v52}{.2}\bow{v42}{v52}{-.2}
\edge{v42}{v53}\bow{v42}{v53}{.2}\bow{v42}{v53}{-.2}
\bow{v43}{v53}{.075}\bow{v43}{v53}{.2}\bow{v43}{v53}{-.075}\bow{v43}{v53}{-.2}
\bow{v43}{v54}{.15}\bow{v43}{v54}{-.15}
\edge{v44}{v54}\bow{v44}{v54}{.075}\bow{v44}{v54}{.2}\bow{v44}{v54}{-.075}\bow{v44}{v54}{-.2}
\edge{v44}{v55}\freetext(10,5){\underline{Level}}\freetext(10,4){0}
\freetext(10,3){1}\freetext(10,2){2}\freetext(10,1){3}\freetext(10,0){4}\freetext(10,-1){5}
\nodetext{v50}(0,-.5){$k=0$}
\nodetext{v51}(0,-.5){$k=1$}\nodetext{v52}(0,-.50){$k=2$}\nodetext{v53}(0,-.50){$k=3$}
\nodetext{v54}(0,-.50){$k=4$}\nodetext{v55}(0,-.50){$k=5$}
\end{graph}
\caption{The Euler graph down to level 5} \label{Graph}
\end{figure}

Define $X$ to be the space of infinite edge paths in this graph:
$X=\{x=(x_n)|n=0,1,2,\dots\}$, each $x_n$ being an edge from level
$n$ to level $n+1$.  The vertex through which $x$ passes at level
$n$ is denoted by $(n,k_n(x))$.  We say that an edge is a \emph{left
turn} if it connects vertices $(n,k)$ and $(n+1,k)$ and a
\emph{right turn} if it connects vertices $(n,k)$ and $(n+1,k+1)$.
$X$ is a metric space in the usual way: $d(x,y)=2^{-j}$, where
$j=\inf\{i|x_i\neq y_i\}$. The cylinder sets, where a finite number
of edges are specified, are clopen sets that generate the topology
of $X$.

\begin{figure}[h]
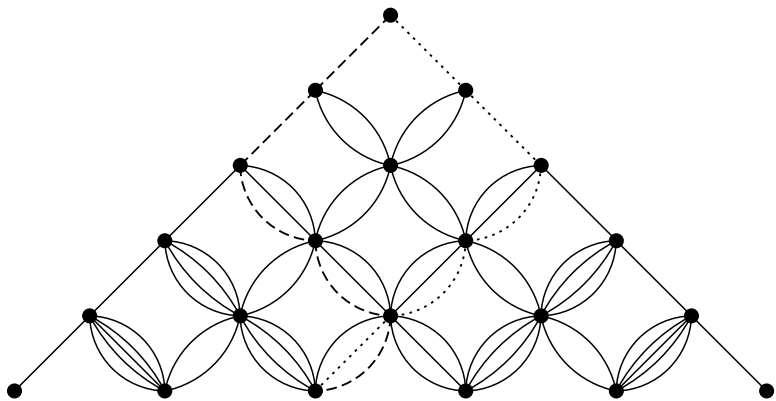

\begin{graph}(7,7)(0,-1.5) \roundnode{v00}(4,4)
\roundnode{v10}(3,3)\roundnode{v11}(5,3)
\roundnode{v20}(2,2)\roundnode{v21}(4,2)\roundnode{v22}(6,2)
\edge{v00}{v10}[\graphlinedash{4
2}\graphlinewidth{.025}]\edge{v00}{v11}[\graphlinedash{1
2}\graphlinewidth{.025}]\edge{v10}{v20}[\graphlinedash{4
2}\graphlinewidth{.025}]\bow{v10}{v21}{.15}
\bow{v10}{v21}{-.15}\bow{v11}{v21}{-.15}
\bow{v11}{v21}{.15}\edge{v11}{v22}[\graphlinedash{1
2}\graphlinewidth{.025}]\roundnode{v30}(1,1)\roundnode{v31}(3,1)\roundnode{v32}(5,1)\roundnode{v33}(7,1)
\edge{v20}{v30}\edge{v20}{v31}\bow{v20}{v31}{.2}\bow{v20}{v31}{-.2}[\graphlinedash{4
2}\graphlinewidth{.025}]\bow{v21}{v31}{.15}\bow{v21}{v31}{-.15}
\bow{v21}{v32}{.15}\bow{v21}{v32}{-.15}\edge{v22}{v32}\bow{v22}{v32}{.2}[\graphlinedash{1
2}\graphlinewidth{.025}]\bow{v22}{v32}{-.2}\edge{v22}{v33}
\roundnode{v40}(0,0)\roundnode{v41}(2,0)\roundnode{v42}(4,0)\roundnode{v43}(6,0)\roundnode{v44}(8,0)
\edge{v30}{v40}\bow{v30}{v41}{.075}\bow{v30}{v41}{.2}\bow{v30}{v41}{-.075}\bow{v30}{v41}{-.2}\bow{v31}{v41}{.15}
\bow{v31}{v41}{-.15}\bow{v31}{v42}{.2}\bow{v31}{v42}{-.2}[\graphlinedash{4
2}\graphlinewidth{.025}]\edge{v31}{v42}\edge{v32}{v42}\bow{v32}{v42}{.2}[\graphlinedash{1
2}\graphlinewidth{.025}] \bow{v32}{v42}{-.2}
\bow{v32}{v43}{.15}\bow{v32}{v43}{-.15}
\bow{v33}{v43}{.075}\bow{v33}{v43}{.2}\bow{v33}{v43}{-.075}\bow{v33}{v43}{-.2}\edge{v33}{v44}
\roundnode{v50}(-1,-1)\roundnode{v51}(1,-1)\roundnode{v52}(3,-1)\roundnode{v53}(5,-1)\roundnode{v54}(7,-1)\roundnode{v55}(9,-1)
\edge{v40}{v50}\edge{v40}{v51}\bow{v40}{v51}{.075}\bow{v40}{v51}{.2}\bow{v40}{v51}{-.075}\bow{v40}{v51}{-.2}
\bow{v41}{v51}{.15}\bow{v41}{v51}{-.15}
\bow{v41}{v52}{.075}\bow{v41}{v52}{.2}\bow{v41}{v52}{-.075}\bow{v41}{v52}{-.2}
\edge{v42}{v52}[\graphlinedash{1
2}\graphlinewidth{.025}]\bow{v42}{v52}{.2}[\graphlinedash{4
2}\graphlinewidth{.025}]\bow{v42}{v52}{-.2}
\edge{v42}{v53}\bow{v42}{v53}{.2}\bow{v42}{v53}{-.2}
\bow{v43}{v53}{.075}\bow{v43}{v53}{.2}\bow{v43}{v53}{-.075}\bow{v43}{v53}{-.2}
\bow{v43}{v54}{.15}\bow{v43}{v54}{-.15}
\edge{v44}{v54}\bow{v44}{v54}{.075}\bow{v44}{v54}{.2}\bow{v44}{v54}{-.075}\bow{v44}{v54}{-.2}
\edge{v44}{v55}
\end{graph}
\caption{The dotted path maps to the dashed path} \label{action}
\end{figure}

  As with other Bratteli diagrams, we define a partial order on
  the edges in the diagram which extends to the entire space of
  paths.  Two edges are comparable if they terminate in the same
  vertex.  For each vertex, totally order the set
  of edges that terminate at that vertex.  This is pictured by agreeing
  that the minimal edge between two vertices is the leftmost edge, and the
  edges increase from left to right.  Two paths, $x,y$, are comparable if
  they agree after some level.  Then $x$ is less than $y$ if the
  last edges of $x$ and $y$ that do not agree, $x_i$ and $y_i$,
  are such that $x_i$ is less than $y_i$ in the edge ordering.
  We define the \emph{adic transformation} $T$ on the set of non-maximal paths in $X$ to map a path to
  the next largest path according to this induced partial order.
  Then two paths are in the same orbit if and only if they are comparable.
Figure \ref{action} shows a particular path $x$ and its image
$y=Tx$.  There are a countable number of paths for which there is no
next largest path, we denote the set of such paths by
$X_{\text{max}}$.  In particular for each $k$ in $\mathbb{N}\cup
\{\infty\}$, we define $x_{\max(k)}$ to be the infinite edge path
that follows the unique finite path from (0,0) to $(k,k)$ and for
$n\geq k$, follows the right most edge connecting vertices $(n,k)$
and $(n+1,k)$.  In a similar fashion there are a countable number of
paths for which there is no next smaller path, denoted by
$X_{\min}$, and for every $k$ in $\mathbb{N}\cup\{\infty\}$ a unique
path $x_{\min(k)}\in X_{\min}$.  The transformation $T$ maps
$x_{\max{k}}$ to $x_{\min(k)}$.

\begin{figure}[h]
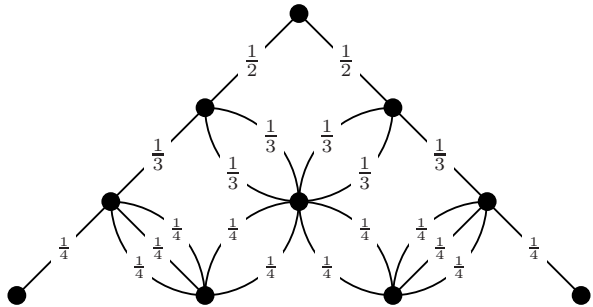


\centering
\begin{graph}(7.5,5)
\unitlength=1.25\unitlength
\roundnode{00}(0,0)\roundnode{20}(2,0)\roundnode{40}(4,0)\roundnode{60}(6,0)
\roundnode{11}(1,1)\roundnode{31}(3,1)\roundnode{51}(5,1)\roundnode{22}(2,2)
\roundnode{42}(4,2)\roundnode{33}(3,3)\edge{33}{00}\edge{33}{60}\bow{22}{31}{-.2}\bow{22}{31}{.2}
\bow{42}{31}{.2}\bow{42}{31}{-.2}\bow{11}{20}{.2}\edge{11}{20}\bow{11}{20}{-.2}\bow{31}{20}{-.2}\bow{31}{20}{.2}
\bow{31}{40}{.2}\bow{31}{40}{-.2}\bow{51}{40}{-.2}\edge{51}{40}\bow{51}{40}{.2}
\edgetext{33}{22}{$\frac{1}{2}$}\edgetext{33}{42}{$\frac{1}{2}$}\edgetext{22}{11}{$\frac{1}{3}$}
\bowtext{22}{31}{.2}{$\frac{1}{3}$}\bowtext{22}{31}{-.2}{$\frac{1}{3}$}\bowtext{42}{31}{.2}{$\frac{1}{3}$}
\bowtext{42}{31}{-.2}{$\frac{1}{3}$}\edgetext{42}{51}{$\frac{1}{3}$}
\edgetext{11}{00}{\tiny $\frac{1}{4}$}\bowtext{11}{20}{-.2}{\tiny
$\frac{1}{4}$}\edgetext{11}{20}{\tiny $\frac{1}{4}$}
\bowtext{11}{20}{.2}{\tiny $\frac{1}{4}$} \bowtext{31}{20}{.2}{\tiny
$\frac{1}{4}$} \bowtext{31}{20}{-.2}{\tiny $\frac{1}{4}$}
\bowtext{31}{40}{.2}{\tiny $\frac{1}{4}$}\bowtext{31}{40}{-.2}{\tiny
$\frac{1}{4}$}\bowtext{51}{40}{.2}{\tiny$\frac{1}{4}$}\edgetext{51}{40}{\tiny$\frac{1}{4}$}
\bowtext{51}{40}{-.2}{\tiny$\frac{1}{4}$}\edgetext{51}{60}{\tiny$\frac{1}{4}$}
\end{graph}
\caption{The symmetric measure given by weights on edges}
\label{SymmetricMeasure}
\end{figure}

  The \emph{symmetric measure}, $\eta$, on the infinite path space $X$ of the Euler graph is the
  Borel probability measure that for each $n$ gives every cylinder of length $n$
  starting at the root vertex
  the same measure.  Clearly $\eta$ is $T$-invariant.  The measure
  of any cylinder set can be computed by multiplying \emph{weights} on the edges,
  each weight on an edge connecting level $n$ to level $n+1$
  being $1/(n+2)$.  We can think of the weights as assigning equal probabilities
  to all the allowed steps for a random walker who starts at the
  root and descends step by step to form an infinite path $x\in
  X$.  The main result of \cite{BKPS}, proved by probabilistic
  methods, is that the symmetric measure is ergodic.

There is a bijective correspondence between paths (or cylinders) of
length $n_0$ starting at the root vertex and terminating at vertex
$(n_0,k_0)$ and permutations of $\{1,2,\dots,n_0+1\}$ with $k_0$
rises. This correspondence is crucial for our proof of Proposition
\ref{Unique}, which is the essential component of our main result.
Consider the cylinder set defined by the single edge connecting the
vertex $(0,0)$ to the vertex $(1,0)$. This cylinder set is of length
1 with 1 left turn, and we assign to it the permutation 21, which
has one fall. Likewise, the cylinder set defined by the single edge
connecting the vertex $(0,0)$ to the vertex $(1,1)$ is of length 1
with one right turn, and we assign to it the permutation 12, which
has one rise.  When a cylinder $F$ of length $n$, corresponding to
the permutation $\pi(F)$ of $\{1,2,\dots, n+1\}$, is extended by an
edge from level $n$ to level $n+1$, we extend $\pi(F)$ in a unique
way to a permutation of $\{1,2,\dots,n+2\}$, as follows.  If $F$ is
extended by a left turn down the $i$'th edge connecting $(n,k)$ to
$(n+1,k)$, insert $n+2$ into $\pi(F)$ in the $i$'th place that adds
an additional fall to the total number of falls of $\pi(F)$.
Likewise, if $F$ is extended by a right turn down the $i$'th edge
connecting $(n,k)$ to $(n+1,k+1)$, insert $n+2$ into $\pi(F)$ in the
$i$'th place that adds an additional rise to the total number of
rises of $\pi(F)$.

\begin{figure}[h]
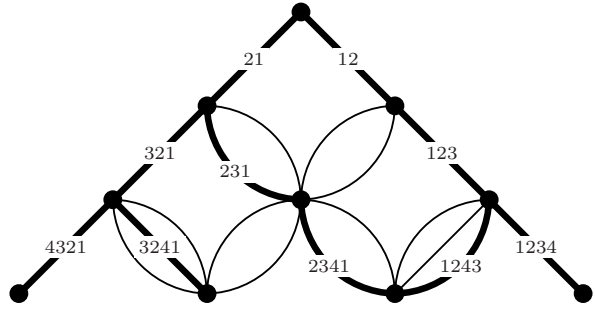
\centering
\begin{graph}(7.5,5)
\unitlength=1.25\unitlength
\roundnode{00}(0,0)\roundnode{20}(2,0)\roundnode{40}(4,0)\roundnode{60}(6,0)
\roundnode{11}(1,1)\roundnode{31}(3,1)\roundnode{51}(5,1)\roundnode{22}(2,2)
\roundnode{42}(4,2)\roundnode{33}(3,3)\edge{33}{00}[\graphlinewidth{.07}]
\edge{33}{60}[\graphlinewidth{.07}]\bow{22}{31}{-.2}[\graphlinewidth{.07}]
\bow{22}{31}{.2}
\bow{42}{31}{.2}\bow{42}{31}{-.2}\bow{11}{20}{.2}\edge{11}{20}[\graphlinewidth{.07}]
\bow{11}{20}{-.2}\bow{31}{20}{-.2}\bow{31}{20}{.2}
\bow{31}{40}{.2}\bow{31}{40}{-.2}[\graphlinewidth{.07}]
\bow{51}{40}{-.2}\edge{51}{40}\bow{51}{40}{.2}[\graphlinewidth{.07}]
\edgetext{33}{22}{\scriptsize 21}\edgetext{33}{42}{\scriptsize
12}\edgetext{22}{11}{\scriptsize 321}\edgetext{42}{51}{\scriptsize
123} \bowtext{22}{31}{-.2}{\scriptsize 231}
\edgetext{11}{00}{\scriptsize 4321}\edgetext{11}{20}{\scriptsize
3241} \bowtext{31}{40}{-.2}{\scriptsize
2341}\bowtext{51}{40}{.2}{\scriptsize
1243}\edgetext{51}{60}{\scriptsize 1234}
\end{graph}
 \caption{Some cylinders and their corresponding permutations}\label{cylinders}
\end{figure}

\begin{example}We shall determine which cylinder $F$ corresponds to the
permutation $\pi(F)=2341$.  When we delete 34, we have the
permutation 21, which means the first edge of $F$ is the unique edge
connecting vertex (0,0) to (1,0).  When the 3 is added, it is added
into the middle of 21, creating the permutation 231, which has one
rise and one fall.  Therefore a rise has been added, and there is no
earlier space that the 3 could have been placed to add a rise to 21.
Therefore, the second edge of $F$ is the first edge connecting
vertex (1,0) to vertex (2,1).  Since $\pi(F)$ has 2 rises and one
fall, the addition of 4 into 2341 adds another rise.  It is added in
the first place it could be added in order to add a rise, meaning
that the third edge of $F$ must be the first edge connecting vertex
(2,1) to (3,2).  See Figure \ref{cylinders}.
\end{example}

This correspondence produces a labeling of infinite paths in the
Euler graph starting at the root; then the path space $X$
corresponds to the set of all linear orderings of
$\mathbb{N}=\{1,2,3,\dots\}$, and the adic transformation $T$ can be
thought of as moving from an ordering to its successor.

\section{Proof of the main result---uniqueness of the symmetric measure}
The \emph{dimension} $\dim(n,k)$ of a vertex $(n,k)$ is defined to
be the number of paths connecting the root vertex $(0,0)$ to the
vertex $(n,k)$. For any cylinder $F$ define $\dim(F,(n,k))$ to be
the number of paths in $F$ that connect the root vertex $(0,0)$ to
the vertex $(n,k)$.  If $F$ is a cylinder starting at the root
vertex $(0,0)$, then  $\dim(F,(n,k))$ is the number of paths from
the terminal vertex of $F$ to $(n,k)$.  Denote by $\mathcal{I}$ the
$\sigma$-algebra of $T$-invariant Borel measurable subsets of $X$.
 The following is a well-known result about the measures of
cylinders in an adic situation.

\begin{lem}[Vershik \cite{Vershik1},\cite{Vershik2}]
Let $\mu$ be an invariant probability measure for the Euler adic
transformation.  Then for every cylinder set $F$ and $\mu$-almost
every $x\in X$,
$$\lim_{n\to
\infty}\frac{\dim(F,(n,k_n(x)))}{\dim((n,k_n(x)))}=E_{\mu}(\chi_F|\mathcal{I}).$$
\label{GeneralErgodicity}
\end{lem}

It is clear from the Markov property and the patterns of weights on
the edges that with respect to $\eta$ almost every path $x\in X$ has
infinitely many left and right turns.
\begin{lem}
Let $\mu$ be an invariant fully-supported ergodic probability
measure for the Euler adic transformation. For $\mu$-almost every
$x\in X$, there are infinitely many left and right turns (i.e.,
$k_n(x)$ and $n-k_n(x)$ are unbounded a.e.).
  \label{MuAlmostEverywhere}
  \end{lem}

  \begin{proof}
  For each $K=1,2,\dots$, let $$A_K=\{x\in X|x \text{ has no more than }K \text{
  right turns}\}.$$  Then
  $A_K$ is a proper closed $T$-invariant set.  Since $\mu$ is ergodic and
  fully supported, $\mu(A_K)=0$.  Similarly, the measure
  of the set of paths with a bounded number of left turns is 0.
   \end{proof}

\begin{prop}
  If $\mu$ is an invariant probability measure for the Euler adic
  transformation
  such that $k_n(x)$ and $n-k_n(x)$ are unbounded a.e., then $\mu=\eta$.
  \label{Unique}
  \end{prop}

\begin{proof}
For any string $w$ on an ordered alphabet denote by $r(w)$ the
number of rises in $w$ and by $f(w)$ the number of falls in $w$
(defined as above).  Let $F$ and $F'$ be cylinder sets in $X$
specified by fixing the first $n_0$ edges, and let $\pi(F)$ and
$\pi(F')$ be the permutations assigned to them by the correspondence
described in the preceding section.  Suppose that the paths
corresponding to $F$ and $F'$ terminate in the vertices $(n_0,k_0)$
and $(n_0,k_0')$ respectively.  Fix $n\gg n_0$ and $k\gg k_0$.

 We will show that $\mu(F)=\mu(F')$, and hence $\mu=\eta$. If
$k_0=k_0'$, then $\mu(F)=\mu(F')$ by the $T$-invariance of $\mu$, so
there is nothing to prove.  We will deal with the case $k_0\neq
k_0'$.

\begin{example}
Let $F$ be the dotted cylinder, and $F'$ the dashed.  Then
$\pi(F)=213$ and $\pi(F')=132$.  For purposes of following an
example through, we will let $(n,k)=(8,4)$.
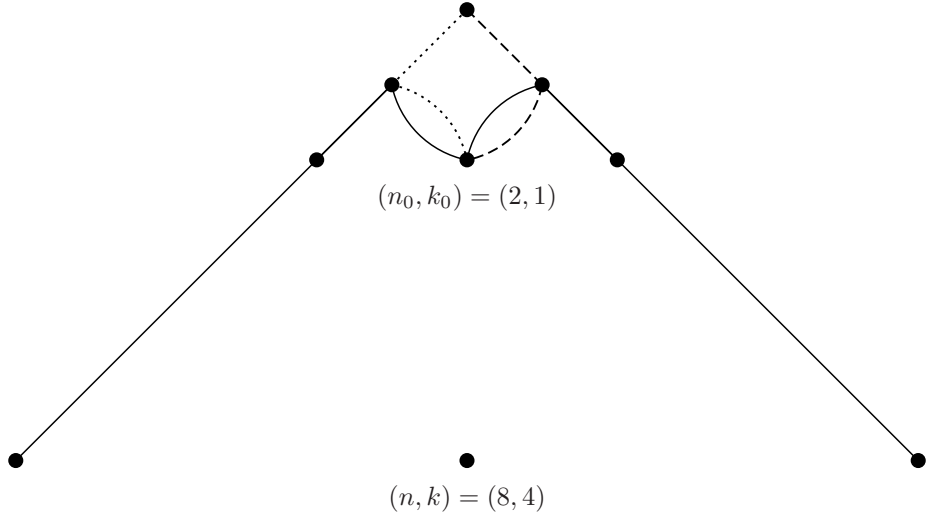
\begin{center}
\begin{graph}(12,7)(-2,-2.5) \roundnode{v00}(4,4)
\roundnode{v10}(3,3)\roundnode{v11}(5,3)
\roundnode{v20}(2,2)\roundnode{v21}(4,2)\roundnode{v22}(6,2)
\edge{v00}{v10}[\graphlinedash{1 2}\graphlinewidth{.025}]
\edge{v00}{v11}[\graphlinedash{4
2}\graphlinewidth{.025}]\edge{v10}{v20}
\bow{v10}{v21}{.15}[\graphlinedash{1 2}\graphlinewidth{.025}]
\bow{v10}{v21}{-.15}\bow{v11}{v21}{-.15} \bow{v11}{v21}{.15}
[\graphlinedash{4 2}\graphlinewidth{.025}]\edge{v11}{v22}
\nodetext{v21}(0,-.5){$(n_0,k_0)=(2,1)$} \roundnode{v50}(-2,-2)
\roundnode{v55}(10,-2)
\edge{v10}{v50}\edge{v11}{v55}\roundnode{v84}(4,-2)\nodetext{v84}(0,-.5){$(n,k)=(8,4)$}
\end{graph}
\end{center}

 \end{example}

  Each path $s$ in $\Gamma$ from $(n_0,k_0)$ to $(n,k)$
  corresponds to a permutation $\sigma_s$ of $\{1,2,\dots,n+1\}$
  with $k$ rises in which $1,2,\dots,n_0+1$ appear in the order
  $\pi(F)$.  Counting $\dim(F,(n,k))$ is equivalent to counting
  the number of distinct such $\sigma_s$.  Each such permutation
  $\sigma_s$ has associated to it a permutation $t(\sigma_s)$ of
  $\{n_0+2,\dots,n+1\}$ obtained by deleting $1,2,\dots,n_0+1$
  from $\sigma_s$, see Example \ref{Ex2}.  Taking a reverse view, one obtains $\sigma_s$
  from $\rho=t(\sigma_s)$ by inserting $1,2,\dots,n_0+1$ from left
  to right, in the order prescribed by $\pi(F)$, into $\rho$.

  We define a \emph{cluster} in $\sigma_s$ to be a subset of
  $\{1,2,\dots,n_0+1\}$ whose members are found consecutively in
  $\sigma_s$, with no elements of $\{n_0+2,\dots,n+1\}$ separating
  them, in the order prescribed by $\pi(F)$.  The set $M_s$ of
  clusters in $\sigma_s$ is an ordered partition of the
  permutation $\pi(F)$, and we define $$r(M_s)=\sum_{c\in
  M_s}r(c).$$  In general, $1\leq |M_s|\leq n_0+1$ and $0\leq
  r(M_s)\leq k_0$.

\begin{example}
  $\pi(F)=$\textbf{213}
  \[\begin{array}{lllll}\sigma_{s_1}=\textbf{2}97\textbf{1}46\textbf{3}85
  &&&&\sigma_{s_2}=96\textbf{2}47\textbf{13}58\\
  t(\sigma_{s_1})=974685 &&&& t(\sigma_{s_2})=964758\\
  M_{s_1}=\{2,1,3\} &&&& M_{s_2}=\{2,13\}\\
  r(M_{s_1})=0 &&&& r(M_{s_2})=1\\
  \end{array}\]\label{Ex2}
  \end{example}

  Given a permutation $\rho$ of $\{n_0+2,\dots,n+1\}$, $0\leq
  m\leq |\rho|+1$, and an ordered partition $M$ of $\pi(F)$ with
  $|M|=m$, there are $C(|\rho|+1,m)$ (the binomial coefficient)
  choices for how to insert the
  members of $M$ as clusters into the permutation $\rho$ in order
  to form a permutation $\sigma_s$.  But not all of these choices
  yield valid permutations $\sigma_s$, which have exactly $k$
  rises.  Looking more closely, we see that placing a cluster
  $c\in M$ at the tail end of $\rho$ or into a rise in $\rho$ produces
  a permutation $\overline{\rho}$ whose number of rises is
  $r(\overline{\rho})=r(\rho)+r(c)$, while placing $c$ at the
  beginning or into a fall produces $\overline{\rho}$ with
  $r(\overline{\rho})=r(\rho)+r(c)+1$.  So we must have
  $$k=r(\sigma_s)=r(\rho)+r(M)+$$ $$\#\{c\in M | c\text
  { is placed into a fall or at the beginning of }\rho\}.$$
In order to count the number of valid ways to place the members of
$M$ into $\rho$, we will first determine how many ways there are to
place the clusters that will create a new rise. There
  are $n-n_0-(k-k_0)+1$ possible places in $\rho$ to place a cluster of
  $M$ and create a new rise
  (the $n-n_0-(k-k_0)$ falls and at the end of $\rho$),
  and since $|\{c\in M | c\text
  { is placed into a fall or at the beginning of }\rho\}|$ is $k-r(\rho)-r(M)$ we must choose
 $l(r(\rho),M)=k-r(\rho)-r(M)$ of them.  For each of these possibilities we must
  then choose places to place the remaining clusters of $M$.  There
  are $k-k_0+1$ places, and $m-(k-r(\rho)-r(M))$
  remaining clusters.  Therefore the number of ways (given $m, M,$ and $r=r(\rho)$) to place the members of $M$ into $\rho$
in such a way  as to form a valid permutation $\sigma_s$, with $k$
rises and $n-k$ falls, is
$$C(n-n_0-(k-k_0)+1,k-r(\rho)-r(M))\,C(k-k_0+1,m-(k-r(\rho)-r(M))).$$

For each $m=1,2,\dots, n_0+1$ denote by $P_m(F)$ the set of ordered
partitions $M$ of $\pi(F)$ such that $|M|=m$.  For each
$r=0,1,\dots,n-n_0-1$ denote by $Q(n,n_0,r)$ the set of permutations
of $\{n_0+2,\dots,n+1\}$ with exactly $r$ rises, so that
$$|Q(n,n_0,r)|=A(n-n_0-1,r).$$

\begin{example} $\pi(F)=213$, $(n,k)=(8,4)$.
\[\begin{array}{lllll}
P_1(F)=\{\{213\}\}\\
P_2(F)=\{\{2,13\},\{21,3\}\}\\
P_3(F)=\{\{2,1,3\}\}\\
Q(8,2,0)=\{987654\} &&&|Q(8,2,0)|=A(5,0)=1\\
Q(8,2,5)=\{456789\} &&&|Q(8,2,5)|=A(5,5)=1\\
\end{array}\]
\end{example}

In order to compute $\dim(F,(n,k))$ we will partition the set of
permutations $\sigma_s$ of $\{1,2,\dots,n+1\}$ which have exactly
$k$ rises and a subpermutation $\pi(F)$ in the following manner.
First, partition the set of all such permutations according to the
cardinality $m\in\{1,\dots,n_0+1\}$ of the set of  clusters.  Now
partition further by grouping the $\sigma_s$ by their corresponding
$M_s\in P_m(F)$. Partition yet again by grouping the $\sigma_s$ by
the number of rises, $r$, in their corresponding permutations
$t(\sigma_s)$ of $\{n_0+2,\dots,n+1\}$. Since the total number of
rises in each $\sigma_s$ must be $k$, and there are $r(M)$ rises in
$M$, the minimum such $r$ is $k-r(M)-m$ and the maximum is $k-r(M)$.
Finally, partition these groups by the distinct permutations
$t(\sigma_s)$ of $\{n_0+2,\dots,n+1\}$ which meet the foregoing
criteria. Then for each $r$ there are $|Q(n,n_0,r)|=A(n-n_0-1,r)$ of
these groups (the same number for each $F, m, M$), and each group
has cardinality
$$C(n-n_0-(k-k_0)+1,k-r-r(M))\,C(k-k_0+1,m-(k-r-r(M))).$$

\begin{example} Fix $\pi(F)=213$ and $(n,k)=(8,4)$.

1. Partition by $m$.
\[\begin{array}{lll}m=1 &:& \sigma_s = \dots 213\dots\\
m=2&:&\sigma_s= \dots 2\dots 13\dots \text{ or }\sigma_s=\dots 21
\dots  3\dots\\
m=3&:&\sigma_s=\dots 2\dots 1\dots 3\dots \end{array}\]

2. Fix $m=2$, partition by $M$.
\[\begin{array}{lll}\dots2\dots 13\dots &\text{ appears in each
}\sigma_s\\
\dots 21\dots 3\dots &\text{ appears in each }\sigma_s\end{array}\]

3.  Fix $m=2$, $M=\{2,13\}$, then $r(M)=1$, partition over
$r(t(\sigma_s))$.
\[\begin{array}{lll}r(t(\sigma_s))=1\\
r(t(\sigma_s))=2\\
r(t(\sigma_s))=3\end{array}\]

4.  Fix $m=2$, $M=\{2,13\}$, $r(M)=1$, $r(t(\sigma_s))=3$, then
partition over $t(\sigma_s)$ in $Q(8,2,3)$.  There are $A(5,3)=302$
sets in this partition, and each one contains $C(4,0)C(4,2)=6$
elements. (To each permutation of 4,5,6,7,8 we can now add only one
rise, and $r(M)$ is already $1$, so all the clusters must be
inserted into rises or at the end.)
\end{example}

For ease of notation, for each $r$, let
$l=l(r,M)=\inf\{|M|,k-r-r(M)\}$ and $p=p(M)=k-r(M)-|M|$. Then
$\dim(F,(n,k))=$
$$\sum_{m=1}^{n_0+1}\sum_{M\in P_m(F)}\sum_{r=p}^{k-r(M)}
A(n-n_0-1,r)C(n-n_0-(k-k_0)+1,l)\,C(k-k_0+1,m-l).$$

Regrouping the terms, rewrite the sum as
\begin{equation}\dim(F,(n,k))=\sum_{r=k-(n_0+1)}^k\left(\sum_{m=1}^{n_0+1}\alpha(F,r,m)\right)A(n-n_0-1,r),
\label{regroup}\end{equation} where we have
$$\alpha(F,r,m)=\sum_{M\in
P_m(F)}C(n-n_0-(k-k_0)+1,l)C(k-k_0+1,m-l).$$

In each term in this sum, the numerator of the first binomial
coefficient factor simplifies to having $l$ factors chosen from
$[n-n_0-(k-k_0)-m,n-n_0-(k-k_0)]$, and the numerator of the second
has $m-l$ factors chosen from $[k-k_0-m,k]$. Combining the two, the
numerator of the binomial factors in each individual term in this
sum consists of $m$ factors chosen from $[n-k-m,n-k]\cup[k-m,k]$,
and the product of the denominators is bounded above by $(m!)^2$.
Since $m\leq n_0+1$, each of the $m$ factors is comparable to (i.e.,
between two constant multiples of) either $k$ or $n-k$.  Then we see
that

$$\beta(F,r)=\sum_{m=1}^{n_0+1}\alpha(F,r,m)$$
is a polynomial in $k$ and $n-k$ of degree $n_0+1$.

Clearly the dominant term of each $\beta(F,r)$ as $k$ and $n-k$ grow
large is the one of maximum degree, $m=n_0+1$. This is exactly the
term $\alpha(F,r,n_0+1)$, which corresponds to the partition $M$ of
$\pi(F)$ into singletons. Then $r(M)=0$, and all the elements of
$\{1,\dots,n_0+1\}$ are placed, as singleton clusters, into a
permutation $\rho$ of $\{n_0+2,\dots,n+1\}$ which has exactly $r$
rises. Since every element of $\{1,\dots,n_0+1\}$ is less than every
element of the set $\{n_0+2,\dots,n\}$, the order of the permutation
$\pi(F)$ has no effect on $\alpha(F,r,n_0+1)$. The \emph{main point}
is that this term is the same for $\pi(F)$ as for any other
permutation $\pi(F')$ of $\{1,2,\dots,n_0+1\}$:
$$\alpha(F,r,n_0+1)=\alpha(F',r,n_0+1)$$
for all $r$ and all $F'$. When $k-r$ of $\{1,\dots,n_0+1\}$ are put
into falls in $\rho$ or at the beginning, and the rest are put into
rises or at the end, {\em no matter which elements are placed in
which slots} we always produce a permutation $\sigma_s$ with exactly
$k$ rises.

Let us now consider the ratio $\dim(F,(n,k))/\dim(F',(n,k))$ when
$n$ and $k$ are both very large.  Divide top and bottom by the sum
on $r$ of the dominant terms (taking maximum degree coefficients in
$k$ and $n-k$ for each $r$),
$$\sum_{r=k-(n_0+1)}^{k}\alpha(F,r,n_0+1)A(n-n_0-1,r),$$
which is the $\emph{same}$ for $F$ and $F'$.  This shows that the
ratio is very near 1 when $k$ and $n-k$ are both very large.

Thus if $k_n(x)$, $n-k_n(x) \to \infty$ a.e. $d\mu$, we have for any
two cylinders $F$ and $F'$ starting at the root vertex and of the
same length that
$$E_{\mu}(X_F|\mathcal{I})(x)=E_{\mu}(X_{F'}|\mathcal{I})(x)\   \text{  a.e. }d\mu.$$
Integrating gives $\mu(F)=\mu(F')$, so that $\mu=\eta$.
\end{proof}

\begin{thm}The symmetric measure $\eta$ is ergodic and is the only fully supported
invariant ergodic Borel probability measure for the Euler adic
transformation.
\end{thm}

\begin{proof} If we show that there is an ergodic measure $\mu$ which has
$k_n(x)$ and $n-k_n(x)$ unbounded a.e., then it will follow from the
Proposition that $\mu=\eta$, and hence $\eta$ is ergodic and is the
only fully-supported $T$-invariant ergodic measure on $X$.

If an ergodic measure has, say, $k_n(x)$ bounded on a set of
positive measure, then $k_n(x)$ is bounded a.e., since each set
$\{x|k_n(x)\leq K\}$ is $T$-invariant.  Let
$\mathcal{E}_0=\emptyset$, and for each $K=1,2,\dots$ let
$\mathcal{E}_K$ be the set of ergodic measures supported on either
$\{x\in X|k_n(x)\leq K\text{ for all }n\}$ or $\{x\in X|n-k_n(x)\leq
K\text{ for all }n\}$.  If no ergodic measure has $k_n(x)$ and
$n-k_n(x)$ unbounded a.e., then the set of ergodic measures is
$$\mathcal{E}=\bigcup_{K}\mathcal{E}_K.$$

Form the ergodic decomposition of $\eta$:
$$\eta=\int_{\mathcal{E}}e\,dP_{\eta}(e)=\sum_{K=1}^{\infty}\int_{\mathcal{E}_K\setminus\mathcal{E}_{K-1}}e\,dP_{\eta}(e).$$

If $S$ is the set of paths $x\in X$ for which both $k_n(x)$ and
$n-k_n(x)$ are unbounded, then, from the remark before Lemma
\ref{MuAlmostEverywhere}, $\eta(S)=1$; but, for each $K$, $e(S)=0$
for all $e$ in $\mathcal{E}_K$. Hence there is an ergodic measure
for which $k_n(x)$ and $n-k_n(x)$ are unbounded a.e..
\end{proof}

\section{Concluding remarks}
\begin{rems}\label{randomperms}{\em  The connection of this Theorem with the statements
made in the Abstract and Introduction about random permutations can
be seen as follows.  As noted above, the space $X$ of infinite paths
in the Euler graph is in correspondence with the set $\mathcal{L}$
of linear orderings of $\mathbb{N}$.  A cylinder set in $X$
determined by fixing an initial path of length $n$ corresponds, as
explained above, to a permutation $\pi_{n+1}$ in the group $S_{n+1}$
of permutations of $\{1,2,\dots , n+1\}$, and thus to the set
$\mathcal{L}(\pi_{n+1})$ of all elements of $\mathcal{L}$ for which
$1,2,\dots , n+1$ appear in the order specified by $\pi_{n+1}$. This
family of clopen cylinder sets defines a compact metrizable topology
and a Borel structure on $\mathcal{L}$.  One way to speak about
``random permutations" would be to give a Borel probability measure
on $\mathcal{L}$.  A Borel probability measure on $X$ is
$T$-invariant if and only if the corresponding measure on
$\mathcal{L}$ assigns, for each $n$ and $0\leq k\leq n$, equal
measure to all cylinders determined by permutations $\pi_{n+1}\in
S_{n+1}$ which have the} same fixed number $k$ of rises.
{\em
According to the Theorem, the only such fully supported measure is the one
determined by the symmetric measure $\eta$ on $X$, which assigns to
each basic set $\mathcal{L}(\pi_{n+1})$, for all $\pi_{n+1}\in
S_{n+1}$,} no matter the number of rises {\em in $\pi_{n+1}$, the
\emph{same} measure, $1/(n+1)!$.}
\end{rems}

\begin{figure}[h]
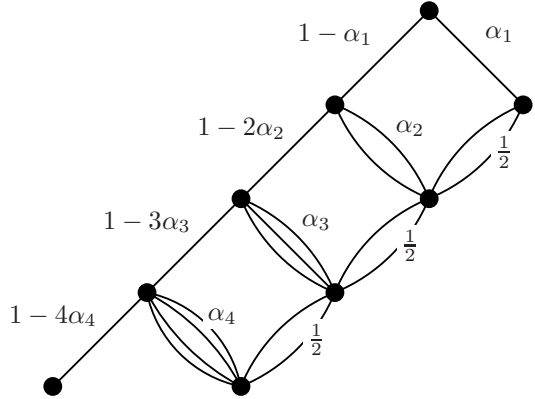

\begin{center}
\begin{graph}(5,5)
\unitlength=1.25\unitlength
\roundnode{1}(4,4)\roundnode{2}(3,3)\roundnode{3}(5,3)\roundnode{4}(2,2)\roundnode{5}(4,2)\roundnode{6}(1,1)
\roundnode{7}(3,1)\roundnode{8}(0,0)\roundnode{9}(2,0)
\edge{1}{2}\freetext(3,3.75){$1-\alpha_1$}\edge{1}{3}\freetext(4.75,3.75){$\alpha_1$}
\edge{2}{4}\freetext(2,2.75){$1-2\alpha_2$}\bow{2}{5}{-.1}\bow{2}{5}{.1}\freetext(3.8,2.75){$\alpha_2$}
\bow{3}{5}{-.1}\bow{3}{5}{.1}\freetext(4.8, 2.5){$\frac{1}{2}$}
\edge{4}{6}\freetext(1,1.75){$1-3\alpha_3$}\bow{4}{7}{-.1}\bow{4}{7}{.1}\edge{4}{7}
\freetext(2.8,1.75){$\alpha_3$}
\bow{5}{7}{-.1}\bow{5}{7}{.1}\freetext(3.8, 1.5){$\frac{1}{2}$}
\edge{6}{8}\freetext(0,.75){$1-4\alpha_4$}\bow{6}{9}{-.075}\bow{6}{9}{.075}\bow{6}{9}{-.15}\bow{6}{9}{.15}
\freetext(1.8,.75){$\alpha_4$}
\bow{7}{9}{-.1}\bow{7}{9}{.1}\freetext(2.8, .5){$\frac{1}{2}$}
\end{graph}
\caption{The invariant measures when $T$ is restricted to the paths
$x$ with vertices $(n,k_n(x))$ with $k_n(x)=0$ or $1$ for all
$n$.}\label{alphas}
\end{center}
\end{figure}

\begin{rems}\label{nonfull}{\em
It is easy to find measures for the Euler adic which are not fully
supported, just by restricting to closed $T$-invariant sets. For
example, if we restrict to the subgraph consisting of paths $x$ with
vertices $(n,k_n(x))$ with $k_n(x)=0$ or $1$ for all $n$, then there
are invariant measures determined by the edge weights in Figure
\ref{alphas}, provided that
$\alpha_{n+1}=\alpha_{n}/(2-2n\alpha_n)$, for example
$\alpha_n=1/2(n+1)$. All such systems are of finite rank (see
\cite{ORW} for the definition). Since measures which are not fully
supported do not involve the full richness of the Euler graph, we
regard them as less interesting than fully supported ergodic
measures.}
\end{rems}

\begin{rems}\label{moreresults}
{\em For dynamic properties of the Euler adic with its symmetric
measure beyond ergodicity,  so far we know that this system is
totally ergodic (there are no roots of unity among its
eigenvalues)\cite{SBF2}, has entropy 0, and is loosely Bernoulli
\cite{SBF2} (see \cite{ORW} for the definition). Determining whether
the system is weakly or even strongly mixing and whether it has
infinite rank are important, challenging, open questions.}
\end{rems}

\begin{figure}[h]
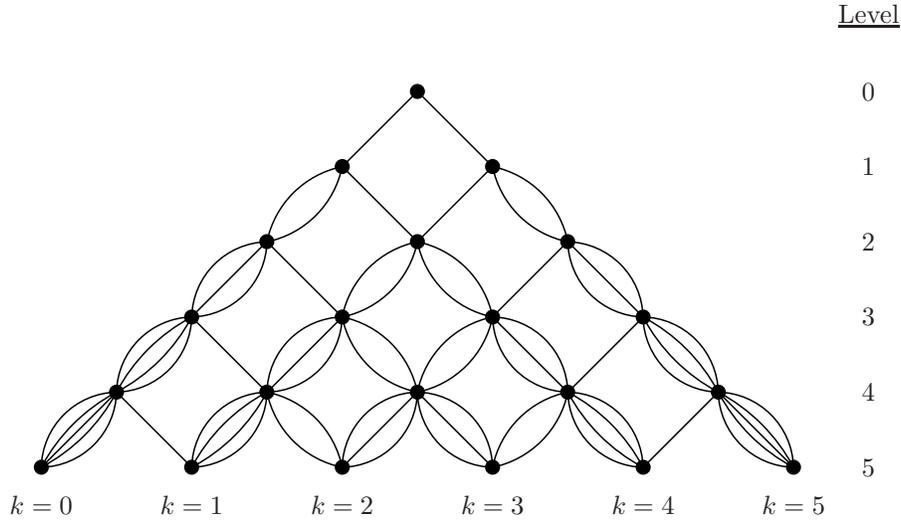

\begin{graph}(9,7)(0,-2) \roundnode{v00}(4,4)
\roundnode{v10}(3,3)\roundnode{v11}(5,3)
\roundnode{v20}(2,2)\roundnode{v21}(4,2)\roundnode{v22}(6,2)\edge{v00}{v10}\edge{v00}{v11}
\edge{v10}{v21}\bow{v10}{v20}{.15}
\bow{v10}{v20}{-.15}\bow{v11}{v22}{-.15}
\bow{v11}{v22}{.15}\edge{v11}{v21}\roundnode{v30}(1,1)\roundnode{v31}(3,1)\roundnode{v32}(5,1)\roundnode{v33}(7,1)
\edge{v20}{v31}\edge{v20}{v30}\bow{v20}{v30}{.2}\bow{v20}{v30}{-.2}\bow{v21}{v32}{.15}\bow{v21}{v32}{-.15}
\bow{v21}{v31}{.15}\bow{v21}{v31}{-.15}\edge{v22}{v33}\bow{v22}{v33}{.2}\bow{v22}{v33}{-.2}\edge{v22}{v32}
\roundnode{v40}(0,0)\roundnode{v41}(2,0)\roundnode{v42}(4,0)\roundnode{v43}(6,0)\roundnode{v44}(8,0)
\edge{v30}{v41}\bow{v30}{v40}{.075}\bow{v30}{v40}{.2}\bow{v30}{v40}{-.075}\bow{v30}{v40}{-.2}\bow{v31}{v42}{.15}
\bow{v31}{v42}{-.15}\bow{v31}{v41}{.2}\bow{v31}{v41}{-.2}\edge{v31}{v41}\edge{v32}{v43}\bow{v32}{v43}{.2}
\bow{v32}{v43}{-.2} \bow{v32}{v42}{.15}\bow{v32}{v42}{-.15}
\bow{v33}{v44}{.075}\bow{v33}{v44}{.2}\bow{v33}{v44}{-.075}\bow{v33}{v44}{-.2}\edge{v33}{v43}
\roundnode{v50}(-1,-1)\roundnode{v51}(1,-1)\roundnode{v52}(3,-1)\roundnode{v53}(5,-1)\roundnode{v54}(7,-1)
\roundnode{v55}(9,-1)
\edge{v40}{v51}\edge{v40}{v50}\bow{v40}{v50}{.075}\bow{v40}{v50}{.2}\bow{v40}{v50}{-.075}\bow{v40}{v50}{-.2}
\bow{v41}{v52}{.15}\bow{v41}{v52}{-.15}
\bow{v41}{v51}{.075}\bow{v41}{v51}{.2}\bow{v41}{v51}{-.075}\bow{v41}{v51}{-.2}
\edge{v42}{v53}\bow{v42}{v53}{.2}\bow{v42}{v53}{-.2}
\edge{v42}{v52}\bow{v42}{v52}{.2}\bow{v42}{v52}{-.2}
\bow{v43}{v54}{.075}\bow{v43}{v54}{.2}\bow{v43}{v54}{-.075}\bow{v43}{v54}{-.2}
\bow{v43}{v53}{.15}\bow{v43}{v53}{-.15}
\edge{v44}{v55}\bow{v44}{v55}{.075}\bow{v44}{v55}{.2}\bow{v44}{v55}{-.075}\bow{v44}{v55}{-.2}
\edge{v44}{v54} \freetext(10,5){\underline{Level}}\freetext(10,4){0}
\freetext(10,3){1}\freetext(10,2){2}\freetext(10,1){3}\freetext(10,0){4}\freetext(10,-1){5}
\nodetext{v50}(0,-.5){$k=0$}
\nodetext{v51}(0,-.5){$k=1$}\nodetext{v52}(0,-.50){$k=2$}\nodetext{v53}(0,-.50){$k=3$}
\nodetext{v54}(0,-.50){$k=4$}\nodetext{v55}(0,-.50){$k=5$}
\end{graph}
\caption{The Reverse Euler graph} \label{reversed}
\end{figure}

\begin{rems}\label{rw}
{\em The Euler adic system  models a negatively reinforced random
walk on the finite graph consisting of two loops based at a single
vertex. In other words, each time an individual loop is followed the
probability that it will be followed at the next stage decreases.
 The adic system based on the Bratteli diagram in
Figure \ref{reversed}, which we call the} reverse Euler system{\em ,
models the positively reinforced random walk (each time an
individual loop is followed the probability that it will be followed
at the next stage increases) on this graph. Identification of the
ergodic measures for such systems and determination of their
dynamical properties has implications for the properties of the
random walks. While the Euler adic has only one fully supported
ergodic measure as we have shown, there is an infinite
(one-parameter) family of fully supported ergodic measures for the
reverse Euler adic. The connection between adic systems of this kind
and reinforced random walks will be explored in a future publication
\cite{FP}.}
\end{rems}

\emph{Acknowledgment}: We thank Ibrahim Salama and Omri Sarig for
useful conversations and the referees for helping us clarify the exposition.








\bibliographystyle{plain}
\bibliography{biblio}

\end{document}